\documentclass[11pt]{article}
\usepackage{setspace,xargs}
\usepackage{amsmath}

\usepackage{graphicx,subfigure}
\usepackage{graphicx,xcolor}
\graphicspath{{./images/}}
 \usepackage{epsfig}

\usepackage [a4paper, bindingoffset=-.5cm, footskip=1cm, headheight = 16pt, top=3cm, bottom=2.5cm, right=3cm, left=3cm ] {geometry}
\input{amssym}
\input{epsf}
\usepackage{stackrel}
\usepackage{algorithm}
\usepackage{algorithmic}

\usepackage{float}
\usepackage{graphicx}
\usepackage{float}
\usepackage{morefloats}
\usepackage{float}
\usepackage{morefloats}
\usepackage{mathtools}

 \usepackage[usenames,dvipsnames]{pstricks}
 \usepackage{epsfig}
 \usepackage{pst-grad} % For gradients
 \usepackage{pst-plot} % For axes
 \usepackage[space]{grffile} % For spaces in paths
 \usepackage{etoolbox} % For spaces in paths
 \makeatletter % For spaces in paths
 \patchcmd\Gread@eps{\@inputcheck#1 }{\@inputcheck"#1"\relax}{}{}
 \makeatother
\usepackage[pagewise]{lineno}
\newtheorem{prealg}{{\bf Algorithm}}

\newtheorem{prethm}{{\bf Theorem}}

\newenvironment{theorem}{\begin{prethm}{\hspace{-0.5
em}{\bf.}}}{\end{prethm}}
\newtheorem{preex}{{\bf Example}}

\newtheorem{prepr}{{\bf Problem}}

\newtheorem{precon}{{\bf Conjecture}}

\newenvironment{conjecture}{\begin{precon}{\hspace{-0.5
em}{\bf.}}}{\end{precon}}
\newtheorem{prelem}{{\bf Lemma}}

\newtheorem{precor}{{\bf Corollary}}

\newtheorem{prepos}{{\bf Proposition}}

\newtheorem{preobv}{{\bf Observation}}

\newtheorem{predef}{{\bf Definition}}

\newtheorem{preproof}{{\bf Proof.}}

\newenvironment{proof}[1]{\begin{preproof}{\rm
               #1}\hfill{$\rule{2mm}{2mm}$}}{\end{preproof}}
               
               \newtheorem{preprooft}{{\bf Proof of Theorem 1.}}

\begin{document}
\renewcommand{\baselinestretch}{1.3}
%\today
\title{{\Large\bf  Star Chromatic Index of Halin Graphs}}

{\small
\author{
%{\sc Marzieh Vahid Dastjerdi}
%footnote{m.vahiddastjerdi@math.iut.ac.ir}
%{\sc Behnaz Omoomi}
%\footnote{bomoomi@cc.iut.ac.ir}
%%\ ,
%and
%%\footnote{e.roshanbin@dal.ca}
%\ and
{\sc  Marzieh Vahid Dastjerdi}}
%\\ [1mm]
%{\small \it Department of Mathematical Sciences}\\
%{\small \it  Isfahan University of Technology} \\
%{\small \it 84156-83111, \ Isfahan, Iran}}
\maketitle
\begin{abstract}
\noindent A star edge coloring of a graph $G$ is a proper edge coloring of $G$ such that  every path and cycle of length four in $G$  uses at least three different colors. The star chromatic index of  $G$, is the smallest integer $k$ for which $G$ admits a star edge coloring with $k$ colors. 
 In this paper,  we obtain tight upper bound $\left\lfloor\frac{3\Delta}{2}\right\rfloor+2$ for the star chromatic index of every Halin graph, that proves the conjecture of Dvo{\v{r}}{\'a}k et al. (J Graph Theory, 72 (2013), 313--326) for cubic Halin graphs. 
\end{abstract}

\section{Introduction}
All graphs considered in this paper are finite and simple. 
For a graph $G$, we denote by $V(G)$ and $E(G)$, its vertex set and edge
set, respectively.
  For every vertex $v\in V(G)$, a vertex $u \in V (G)$ is a {\it neighbor} of $v$ or {\it adjacent} to $v$ if $uv \in E(G)$.
We denote the set of neighbors of $v$ in $G$, by $N_G(v)$.
 Two edges are said to be {\it adjacent} if they are incident to a common vertex. 
 A {\it proper edge coloring} of $G$ is a function assigning  colors to the edges of $G$  in such a way that no two adjacent edges receive the same color. 
 A {\it star edge coloring} of $G$ is a proper edge coloring of $G$ with no bicolored path nor cycle on four edges. 
 The star chromatic index of $G$, denoted by $\chi^\prime_s(G)$, is the minimum $k$ such that $G$ admits a star edeg coloring with $k$ colors.
%A subgraph $F$ of $G$ is said to be \textit{bi-colored}, if the restriction of the proper vertex (edge) coloring of $G$ to $F$, is a vertex (edge) coloring with  at most two colors. 
%Under additional constraints on the proper vertex (edge) coloring of graphs, we get a variety of colorings such as the star vertex  and the star edge coloring. A \textit{star vertex coloring} of  $G$, is a proper vertex coloring such that no path or cycle  on four vertices in G is bi-colored (uses at most two colors) \rm{\cite{coleman,fertin}}.  
%\begin{definition} \rm{\cite{mohar,delta}} A \textit{star edge coloring}   of  a graph $G$ is a proper edge coloring of $G$ such that  no path or cycle of length four in $G$ is bi-colored.    The smallest integer $k$   that $G$ admits a $k$-star edge coloring is called \textit{star chromatic index}  of $G$ and is denoted by  $\chi^\prime_s(G)$.
%\end{definition}

 In 2008, Liu and Deng \cite{delta} introduced the concept of star edge coloring  motivated by the vertex version \rm{\cite{coleman,fertin}}.
  Liu and Deng \cite{delta} presented   upper bound $\left\lceil16(\Delta-1)^\frac{3}{2}\right\rceil$ on the star chromatic index of graphs with maximum  degree $\Delta\geq 7$.
 In \cite{mohar},  Dvo{\v{r}}{\'a}k et al. obtained a near-linear upper bound $\Delta. 2^{O(1)\sqrt{\log(\Delta)}}$ for graphs with maximum degree $\Delta$ and provided some upper and lower bounds for complete graphs (see also \cite{Kerdjoudj}). They also considered
cubic graphs and showed that the star chromatic index of such graphs lies between
4 and 7.  Dvo{\v{r}}{\'a}k et al. perposed the following conjecture for subcubic graphs (graphs with $\Delta\leq 3$) that is also opened. For more results on the star chromatic index of subcubic graph see \cite{ Lie,luzar, pradeep}.
\begin{conjecture} \rm{\cite{mohar}}\label{con:cubic}
If $G$ is a subcubic graph, then $\chi^\prime_s(G)\leq 6$. 
\end{conjecture}

 In \rm{\cite{class}}, Bezegov{\'a}  et al.   obtained some upper bounds on the star chromatic index of subcubic outerplanar graphs, trees and outerplanar graphs, as follows.
 \begin{theorem} \rm{\cite{mohar}}
 Let $G$ be a graph with maximum degree $\Delta$.
 \begin{itemize}
 \item[\rm{(}$i$\rm{)}]
 If $G$ is  tree, then $\chi^\prime_s(G)\leq \left\lfloor\frac{3\Delta}{2}\right\rfloor$.
 \item[\rm{(}$ii$\rm{)}]
 If $G$ is  subcubic outerplanar, then $\chi^\prime_s(G)\leq 5$.
 \item[\rm{(}$iii$\rm{)}]
 If $G$ is outerplanar, then $\chi^\prime_s(G)\leq \left\lfloor\frac{3\Delta}{2}\right\rfloor+12$.
 \end{itemize}
 \end{theorem}
 They also conjectured that the star chromatic index of every outerplanar graph is at most $ \left\lfloor\frac{3\Delta}{2}\right\rfloor+1$.
 Wang et al. in \cite{wang_partition,wang_delta4} proved bounds for graphs $G$ under different planarity
conditions and when maximum degree $G$ is at most four. 
\begin{theorem} \rm{\cite{wang_partition,wang_delta4}}\\
$\bullet$ Let G be a planar graph with girth $g$ (length of shortest cycle).
\begin{enumerate}
\item[\rm{(}$i$\rm{)}]
$\chi^\prime_s(G)\leq \frac{11}{4}\Delta+18$.
\item[\rm{(}$ii$\rm{)}]
If $G$ is $K_4$-minor free, then $\chi^\prime_s(G)\leq \frac{9}{4}\Delta+6$.
\item[\rm{(}$iii$\rm{)}]
If $G$ has no 4-cycle, then $\chi^\prime_s(G)\leq \left\lfloor\frac{3\Delta}{2}\right\rfloor+18$.
\item[\rm{(}$iv$\rm{)}]
If $g\geq 5$, then $\chi^\prime_s(G)\leq \left\lfloor\frac{3\Delta}{2}\right\rfloor+13$.
\item[\rm{(}$v$\rm{)}]
If $g\geq 8$, then $\chi^\prime_s(G)\leq \left\lfloor\frac{3\Delta}{2}\right\rfloor+3$.
\item[\rm{(}$vi$\rm{)}]
If $G$ is outerplanar, then $\chi^\prime_s(G)\leq \left\lfloor\frac{3\Delta}{2}\right\rfloor+5$.
\end{enumerate}
\noindent$\bullet$ Let $G$ be a graph with maximum degree $\Delta\leq 4$.
\begin{enumerate}
\item[\rm{(}$i$\rm{)}]
$\chi^\prime_s(G)\leq 14$.
\item[\rm{(}$ii$\rm{)}]
If $G$ is bipartite, then $\chi^\prime_s(G)\leq 13$.
\end{enumerate}
\end{theorem}

\par In this paper, we investigate the star edge coloring of Halin graphs as a family of planar graphs and find upper bound $\left\lfloor\frac{3\Delta}{2}\right\rfloor+2$ for the star chromatic index of this family of planar graphs. This upper bound proves Conjecture~\ref{con:cubic} for the cubic Halin graphs.
 A {\it Halin graph} $G$  is a plane graph consisting of a tree $T$ and  cycle $C$ that each vertex of $T$ is either of degree 1, called {\it leaf}, or of degree at least 3, and $C$ connecting the leaves of $T$ such that $C$ is the boundary of the exterior face. 
The tree $T$ and the cycle $C$ are called the {\it characteristic tree} and the {\it adjoint cycle} of $G$, respectively. 
We usually write  $G=T\cup C$ to make the characteristic tree and the adjoint cycle of $G$ explicit.

In this paper, the most important point in a star edge coloring of a Halin graph is the coloring of its characteristic tree. For this purpose, we consider the characteristic tree as a rooted tree, and in a specific order on its edges, we present the star edge coloring of the tree.
A {\it rooted tree} is a tree in which one vertex has been designated as the {\it root}.
%and every edge is directed away from the root. 
Suppose that $T$ is a  rooted tree. 
If $v$ is a vertex in $T$ other than the root, the {\it parent} of $v$ is the unique vertex that is the eventual predecessor of  $v$, regarding the root. We use $p(v)$ to denote the parent of $v$.

%A {\it rooted tree} is a tree in which one vertex has been designated as the root.
% The {\it height} of a rooted tree $T$ is the number of edges on the longest path between the root and a leaf.
%  The {\it level} of a vertex in $T$ is the distance between the vertex and the root plus one.
%   Note that the level of the root is one. The \{it parent} of a vertex in level $\ell>1$ is its neighbour in level $\ell-1$.  
%   Let $\sigma$ is an ordering of the neighbours of vertex $v$ in $T$. We denote the index of each neighbours $u$ of $v$ in $\sigma$ with $I_{\sigma}(u)$.
%\par This paper is organized as follows. 
%%%%%%%%%%%%%%%%%%%%%%%%%%%%%%%%%%%%%%%%%%%%%%%%%%%%%%%
%%%%%%%%%%%%%%%%%%%%%%%%%%%%%%%%%%%%%%%%%%%%%%%%%%%%%
\section{Star edge coloring of cubic Halin graphs}
In this section, we consider cubic Halin graphs and find tight upper bound 6 for the star chromatic index of these graph. The following theorem proves    Conjecture~\ref{con:cubic} for cubic Halin graphs.
\begin{theorem}\label{Halin-cubic}
If $G=T\cup C$ is a cubic Halin graph, then $\chi^\prime_s(G)\leq 6$.
\end{theorem}
\begin{proof}
{Let $G=T\cup C$ be a cubic Halin graph and  $P : u_0,\ldots,u_\ell$ be a longest path in $T$. It is easy to see that
\[\chi^\prime_s(G)=\begin{cases}
5&~~\text{if}~\ell=2,\\
6&~~\text{if}~\ell=3.
\end{cases}
\]
Therefore, we assume that $P$ is of length at least four and prove the theorem by induction on the length of $C$.  %G% with $\chi^\prime with $\ell>3$. 
By our assumption that $P$ is a longest path, all neighbors of $u_1$, except $u_2$, must be leaves. 
We may change notations to let $v=u_1$, $u:=u_2$,  $w:=u_3$, and $v_1, v_2$, be  the neighbors of $v$ on $C$ as demonstred in Figure~\ref{fig:halin1}.
Let  $x_1, x_2, y_1, y_2$ be the vertices on $C$, where $x_1$ is adjacent to $v_1$ and $x_2$; $y_1$ is adjacent to $v_2$ and $y_2$. 
 Let $x_3$ and $y_3$ be vertices not on $C$, where $x_1x_3$ and $y_1y_3$ are edges in $T$.
 
 Since $u$ is a vertex of degree 3, there exists a path $P^\prime$ in $T$ from $u$ to $x_1$ or from $u$ to $y_1$ with $P\cap P^\prime = \{u\}$.
  Without loss of generality, we  assume that $P^\prime$ is from $u$ to $y_1$. 
Since $P$ is a longest path,  $|P^\prime|\leq2$. 
  Thus, either $u = y_3$ or $u$ is adjacent to $y_3$. 
  
 Let $G^\prime$ be the graph obtained from $G$ by deleting $v,v_1, v_2, y_1$ and adding two new edges $ux_1$ and $uz$, where $z\in V(C)$ is an adjacent  to $y_1$, and $z\neq v_2$. 
 Since $\ell \geq 4$, $G^\prime$ is a Halin graph. We write $G^\prime = T ^\prime \cup C^\prime$.
 
By the inductive hypothesis as $|C^\prime| < |C|$, $\chi^\prime_s(G^\prime)\leq 6$. 
Let $\mathcal{C}=\{1,\ldots,6\}$, and  $\phi^\prime$ is a star edge coloring of $G^\prime$ with color set $\mathcal{C}$.
 Without loss of generality, assume that $\phi^\prime(ux_1)=1$, $\phi^\prime(uz)=2$, and $\phi^\prime(uw)=3$.
 We obtain a 6-star edge coloring $\phi$ of $G$ as follows. 

For every edge that belongs to $E(G)\cap E(G^\prime)$, we set $\phi(e)=\phi(e^\prime)$.
For edges in  $E(G) \setminus E(G^\prime)$, we consider different cases, and in each case, we give  star edge coloring $\phi$ for
$G$ with at most 6 colors. Suppose that $\phi^\prime(w)=\{3,r_1,r_2\}$, $\phi^\prime(x_1)=\{s_1,s_2\}$, and $\phi^\prime(z)=\{t_1,t_2\}$.

\noindent{\bf Case 1:} $u=y_3$.
\\
In this case, we set $\phi(v_1x)=1$, $\phi(v_1v_2)=3$, $\phi(vv_1)=2$, $\phi(uy_1)=1$, and $\phi(uz)=2$.  Let $a\in \mathcal{C}\setminus\{1,2,3,r_1,r_2\}$ and $b\in \mathcal{C}\setminus\{1,2,3,t_1,t_2\}$. We set $\phi(uv)=a$, $\phi(v_2y_1)=b$, and choose an arbitrary color in  $\mathcal{C}\setminus\{1,2,3,a,b\}$ for edge $vv_2$ (see Case~1 in Figure~\ref{fig:halin1}).   
 
\noindent{\bf Case 2:}  $u\neq y_3$. 
\\
We set $\phi(uy_3)=1$, $\phi(uv)=2$, $\phi(v_1x_1)=1$, $\phi(y_1y_2)=3$, and $\phi(y_2z)=2$. Then, we set $\phi(y_2y_3)=a$, $\phi(vv_1)=b$, $\phi(v_2y_1)=c$,  and $\phi(y_3y_1)=e$, and $\phi(vv_2)=d$,  where 
\begin{align*}
a&\in \mathcal{C}\setminus\{1,2,3,t_1,t_2\},~~~~~b\in \mathcal{C}\setminus\{1,2,3,s_1,s_2\},~~~~~
c\in \mathcal{C}\setminus\{1,2,3,a,b\},\\
d&\in \mathcal{C}\setminus\{1,2,3,b,c\},~~~~~
e\in \mathcal{C}\setminus\{1,2,3,a,c\}.
\end{align*}
Note that this colors are choosen in ordering $a,b,c,d,e$.
}
\end{proof}

\vspace*{1cm}
\begin{figure}[!ht]
\begin{center}
\psscalebox{0.8 0.8} % Change this value to rescale the drawing.
{
\begin{pspicture}(0,-3.3947747)(14.809677,3.3947747)
\definecolor{colour0}{rgb}{1.0,0.8,0.8}
\definecolor{colour1}{rgb}{0.8,0.8,0.8}
\psline[linecolor=black, linewidth=0.04, linestyle=dashed, dash=0.17638889cm 0.10583334cm](12.609677,0.3847126)(12.4096775,-1.9152874)
\psline[linecolor=black, linewidth=0.04](5.094526,0.6544096)(5.8399806,1.239258)
\psbezier[linecolor=black, linewidth=0.04](5.3999214,3.3749566)(-3.0147128,2.2305663)(0.35114083,-2.3469946)(5.3999214,-2.346994702671375)
\psline[linecolor=black, linewidth=0.04](2.5389457,0.5139809)(6.0096774,0.5139809)
\psline[linecolor=black, linewidth=0.04](2.5389457,0.5139809)(0.60967743,1.1847126)
\psline[linecolor=black, linewidth=0.04](2.370653,0.5139809)(0.5194335,-0.15918982)
\psline[linecolor=black, linewidth=0.04](1.6072384,2.0017858)(2.8096774,1.8847126)
\psline[linecolor=black, linewidth=0.04](3.8096774,0.4847126)(1.5096774,-1.1152874)
\pscircle[linecolor=black, linewidth=0.04, fillstyle=solid,fillcolor=colour0, dimen=outer](3.0694335,1.8444687){0.3}
\pscircle[linecolor=black, linewidth=0.04, fillstyle=solid,fillcolor=colour0, dimen=outer](2.4547994,0.59812725){0.3}
\pscircle[linecolor=black, linewidth=0.04, fillstyle=solid,fillcolor=colour0, dimen=outer](3.8011408,0.59812725){0.3}
\pscircle[linecolor=black, linewidth=0.04, fillstyle=solid,fillcolor=colour0, dimen=outer](4.8974824,0.59812725){0.3}
\rput[bl](2.3023603,0.51885897){$v$}
\rput[bl](3.648702,0.51885897){$u$}
\rput[bl](4.763336,0.51885897){$w$}
\pscircle[linecolor=black, linewidth=0.04, fillstyle=solid,fillcolor=colour1, dimen=outer](0.43528718,1.2712979){0.3}
\pscircle[linecolor=black, linewidth=0.04, fillstyle=solid,fillcolor=colour1, dimen=outer](1.4450433,2.1127615){0.3}
\pscircle[linecolor=black, linewidth=0.04, fillstyle=solid,fillcolor=colour1, dimen=outer](1.4450433,-1.2530923){0.3}
\pscircle[linecolor=black, linewidth=0.04, fillstyle=solid,fillcolor=colour1, dimen=outer](2.623092,-1.7579703){0.3}
\pscircle[linecolor=black, linewidth=0.04, fillstyle=solid,fillcolor=colour1, dimen=outer](2.623092,2.785932){0.3}
\pscircle[linecolor=black, linewidth=0.04, fillstyle=solid,fillcolor=colour1, dimen=outer](0.43528718,-0.24333617){0.3}
\rput[l](0.33442172,1.2898661){$v_1$}
\rput[bl](0.20967741,-0.3957752){$v_2$}
\rput[bl](1.298387,2.0363255){$x_1$}
\rput[bl](2.4340677,2.7066638){$x_2$}
\rput[bl](1.2877262,-1.3665069){$y_1$}
\rput[bl](2.4877262,-1.8665069){$z$}
\rput[bl](2.9340677,1.7066638){$x_3$}
\pscircle[linecolor=black, linewidth=0.04, fillstyle=solid,fillcolor=colour1, dimen=outer](3.6011407,-2.0774825){0.3}
\rput[l](5.4344215,0.28986603){$r_2$}
\rput[l](5.2344217,1.089866){$r_1$}
\psline[linecolor=black, linewidth=0.04, linestyle=dashed, dash=0.17638889cm 0.10583334cm](3.7096775,0.7847126)(1.583871,1.9137449)
\psline[linecolor=black, linewidth=0.04, linestyle=dashed, dash=0.17638889cm 0.10583334cm](3.7225807,0.3653578)(2.683871,-1.489481)
\rput[bl](0.6,1.7330997){1}
\rput[bl](2.2,-0.4669003){1}
\rput[bl](1.9,2.6330998){$s_1$}
\rput[bl](2.3,2.0330997){$s_2$}
\rput[bl](2.9,-2.3669002){$t_1$}
\rput[bl](0.0,0.33309972){3}
\rput[bl](1.9,-1.8669003){2}
\psbezier[linecolor=black, linewidth=0.04](14.199922,3.3749566)(5.7852874,2.2305663)(9.151141,-2.3469946)(14.199922,-2.346994702671375)
\psline[linecolor=black, linewidth=0.04](11.338945,0.5139809)(14.809677,0.5139809)
\psline[linecolor=black, linewidth=0.04](11.338945,0.5139809)(9.4096775,1.1847126)
\psline[linecolor=black, linewidth=0.04](11.170653,0.5139809)(9.319433,-0.15918982)
\psline[linecolor=black, linewidth=0.04](10.407238,2.0017858)(11.609677,1.8847126)
\psline[linecolor=black, linewidth=0.04](10.407238,-1.173824)(11.985287,-0.8006532)
\psline[linecolor=black, linewidth=0.04](12.609677,0.4847126)(12.009678,-0.61528736)
\psline[linecolor=black, linewidth=0.04](11.885287,-0.90065324)(11.509678,-1.6152874)
\pscircle[linecolor=black, linewidth=0.04, fillstyle=solid,fillcolor=colour0, dimen=outer](11.869433,-0.8165069){0.3}
\pscircle[linecolor=black, linewidth=0.04, fillstyle=solid,fillcolor=colour0, dimen=outer](11.869433,1.8444687){0.3}
\pscircle[linecolor=black, linewidth=0.04, fillstyle=solid,fillcolor=colour0, dimen=outer](11.2548,0.59812725){0.3}
\pscircle[linecolor=black, linewidth=0.04, fillstyle=solid,fillcolor=colour0, dimen=outer](12.601141,0.59812725){0.3}
\pscircle[linecolor=black, linewidth=0.04, fillstyle=solid,fillcolor=colour0, dimen=outer](13.747482,0.59812725){0.3}
\rput[bl](11.20236,0.51885897){$v$}
\rput[bl](12.448702,0.51885897){$u$}
\rput[bl](13.563336,0.51885897){$w$}
\pscircle[linecolor=black, linewidth=0.04, fillstyle=solid,fillcolor=colour1, dimen=outer](9.235287,1.2712979){0.3}
\pscircle[linecolor=black, linewidth=0.04, fillstyle=solid,fillcolor=colour1, dimen=outer](10.245043,2.1127615){0.3}
\pscircle[linecolor=black, linewidth=0.04, fillstyle=solid,fillcolor=colour1, dimen=outer](10.245043,-1.2530923){0.3}
\pscircle[linecolor=black, linewidth=0.04, fillstyle=solid,fillcolor=colour1, dimen=outer](11.423092,-1.7579703){0.3}
\pscircle[linecolor=black, linewidth=0.04, fillstyle=solid,fillcolor=colour1, dimen=outer](11.423092,2.785932){0.3}
\pscircle[linecolor=black, linewidth=0.04, fillstyle=solid,fillcolor=colour1, dimen=outer](9.235287,-0.24333617){0.3}
\rput[l](9.134421,1.2898661){$v_1$}
\rput[bl](9.109677,-0.3957752){$v_2$}
\rput[bl](10.098387,2.0363255){$x_1$}
\rput[bl](11.234068,2.7066638){$x_2$}
\rput[bl](10.087727,-1.3665069){$y_1$}
\rput[bl](11.287726,-1.8665069){$y_2$}
\rput[bl](11.687726,-0.9665069){$y_3$}
\rput[bl](11.634068,1.7066638){$x_3$}
\pscircle[linecolor=black, linewidth=0.04, fillstyle=solid,fillcolor=colour1, dimen=outer](12.401141,-2.0774825){0.3}
\rput[l](13.134421,0.68986607){3}
\psline[linecolor=black, linewidth=0.04, linestyle=dashed, dash=0.17638889cm 0.10583334cm](12.4096775,0.7847126)(10.383871,1.9137449)
\rput[bl](9.4,1.7330997){1}
\rput[bl](11.7,-2.2669003){2}
\rput[bl](10.7,2.6330998){$s_1$}
\rput[bl](11.1,2.0330997){$s_2$}
\rput[bl](10.5,-1.8669003){3}
\rput[bl](11.0,-0.86690027){$e$}
\rput[bl](11.7,-1.5669003){$a$}
\rput[bl](11.9,1.2330997){1}
\rput[bl](11.8,0.63309973){2}
\psline[linecolor=black, linewidth=0.04](3.8096774,-1.2152873)(2.8096774,-1.6152874)
\pscircle[linecolor=black, linewidth=0.04, fillstyle=solid,fillcolor=colour0, dimen=outer](3.8011408,-1.2018727){0.3}
\rput[bl](3.1,-1.3669003){$t_2$}
\rput[bl](1.3,1.0330998){2}
\rput[bl](3.1,1.2330997){1}
\rput[bl](3.4,-0.5669003){2}
\rput[bl](4.3,0.63309973){3}
\rput[bl](1.0,0.23309971){$c$}
\rput[bl](0.6,-0.9669003){$b$}
\rput[bl](3.0,0.63309973){$a$}
\rput[bl](12.281472,-2.1614413){$z$}
\psline[linecolor=black, linewidth=0.04](13.894526,0.7544096)(14.63998,1.3392581)
\rput[l](13.934422,1.1898661){$r_1$}
\rput[l](14.234422,0.28986603){$r_2$}
\psline[linecolor=black, linewidth=0.04](13.609677,-1.5152874)(12.609677,-1.9152874)
\pscircle[linecolor=black, linewidth=0.04, fillstyle=solid,fillcolor=colour0, dimen=outer](13.601141,-1.5018728){0.3}
\pscircle[linecolor=black, linewidth=0.04, fillstyle=solid,fillcolor=colour1, dimen=outer](13.601141,-2.2774825){0.3}
\rput[bl](12.8,-1.6669003){$t_2$}
\rput[bl](12.8,-2.5669003){$t_1$}
\rput[bl](8.8,0.43309972){3}
\rput[bl](12.0,-0.2669003){1}
\rput[bl](12.7,-0.4669003){2}
\rput[bl](9.5,-1.0669003){$c$}
\rput[bl](10.2,-0.06690029){$d$}
\rput[bl](10.1,1.0330998){$b$}
\rput[bl](2.8096774,-3.3947747){\text{Case1}}
\rput[bl](11.709678,-3.3947747){\text{Case 2}}
\end{pspicture}
}
\caption{ The neighborhood around one end of the longest path $P$.}
\label{fig:halin1}
\end{center}
\end{figure}
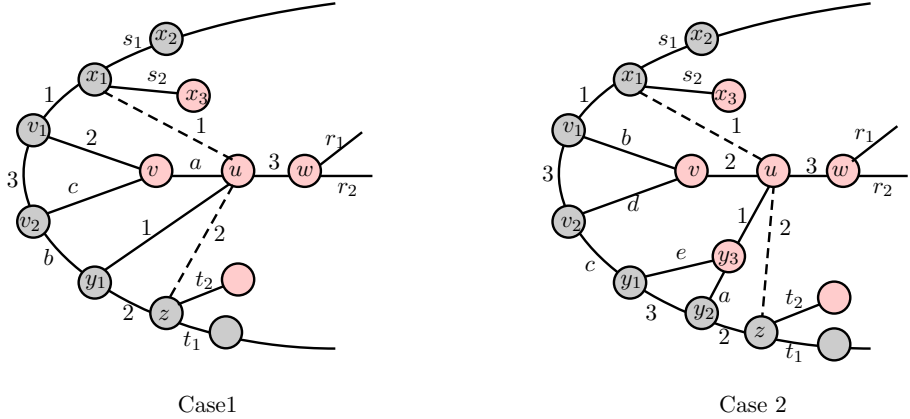

%%%%%%%%%%%%%%%%%%%%%%%%%%%%%%%%%%%%%%%%%%%%%%%%%%%%%%%
%%%%%%%%%%%%%%%%%%%%%%%%%%%%%%%%%%%%%%%%%%%%%%%%%%%%%%%%
\section{Star edge coloring of Halin graphs with $\Delta\geq 4$}
In every Halin graph $G=T\cup C$, we partite the edges of cycle $C$ in to two subsets, as follow.
\begin{align*}
C_s&=\{uv\in E(C):  u~\text{and}~v~\text{have the same neighbor in}~T\},\\
C_d&= E(C)\setminus C_s
\end{align*}
 The first set contains all edges of $C$  connecting two leaves of $T$ with the same neighbor, and the second one  with  different neighbors.
 
 For every non-leaf vertex $v$ of $T$, we define $LN_d(v)$ as the subset of all leaf neighbors of $v$ that  at least an edge of $C_d$ incident to each of them. 
% For each non-leaf vertex $v$ of $T$, we define $V_d(v)$  the subset of $V(C)$ which each vertex in $V_d(v)$ is adjacent to $v$ and also is an end of at least one edges in $C_d$. 
It is easy to see that $|LN_d(v)|\leq \left\lfloor\frac{\Delta}{2}\right\rfloor+1$ (see  Figure~\ref{fig2:LND}).
%, we demostrate two halin graphs with maximum degrees $5$ and $6$ that  the size of  $LN_d(v)$ is maximum in both cases.  

\begin{figure}[H]
\begin{center}
\psscalebox{1.0 1.0} % Change this value to rescale the drawing.
{
\begin{pspicture}(0,-2.2933333)(14.316,2.2933333)
\definecolor{colour0}{rgb}{0.2,0.0,0.8}
\psline[linecolor=black, linewidth=0.04](11.062667,0.5466667)(10.142667,1.9066666)
\psline[linecolor=black, linewidth=0.04](11.142667,0.46666667)(12.102667,1.9466667)
\psline[linecolor=black, linewidth=0.04](0.636,-0.42666668)(1.356,0.53333336)
\psellipse[linecolor=black, linewidth=0.04, dimen=outer](3.18,0.47333333)(3.024,1.62)
\psline[linecolor=black, linewidth=0.04](0.156,0.47333333)(6.204,0.47333333)
\psline[linecolor=black, linewidth=0.04](3.18,0.47333333)(2.1,-1.0386667)
\psline[linecolor=black, linewidth=0.04](3.18,0.47333333)(4.044,-1.0386667)
\psline[linecolor=black, linewidth=0.04](3.18,0.47333333)(3.18,2.0933332)
\psframe[linecolor=black, linewidth=0.04, fillstyle=solid,fillcolor=colour0, dimen=outer](2.208,-0.9306667)(1.872,-1.2666667)
\psframe[linecolor=black, linewidth=0.04, fillstyle=solid,fillcolor=colour0, dimen=outer](4.272,-0.9306667)(3.936,-1.2666667)
\psframe[linecolor=black, linewidth=0.04, fillstyle=solid,fillcolor=colour0, dimen=outer](3.336,2.2933333)(2.976,1.9333333)
\pscircle[linecolor=black, linewidth=0.04, fillstyle=solid,fillcolor=red, dimen=outer](3.15,0.44333333){0.138}
\rput{-8.973528}(-0.008711706,0.77564603){\pscircle[linecolor=black, linewidth=0.04, fillstyle=solid,fillcolor=red, dimen=outer](4.938,0.44333333){0.138}}
\pscircle[linecolor=black, linewidth=0.04, fillstyle=solid,fillcolor=red, dimen=outer](1.314,0.50333333){0.138}
\rput[bl](3.036,0.053333335){$v$}
\pscircle[linecolor=black, linewidth=0.04, fillstyle=solid,fillcolor=red, dimen=outer](6.198,0.44333333){0.138}
\rput{-3.0637553}(-0.023497652,0.00800937){\pscircle[linecolor=black, linewidth=0.04, fillstyle=solid,fillcolor=red, dimen=outer](0.138,0.44333333){0.138}}
\pscircle[linecolor=black, linewidth=0.04, fillstyle=solid,fillcolor=red, dimen=outer](0.618,-0.39666668){0.138}
\psline[linecolor=black, linewidth=0.04](5.436,1.6133333)(4.956,0.53333336)
\pscircle[linecolor=black, linewidth=0.04, fillstyle=solid,fillcolor=red, dimen=outer](5.418,1.5233333){0.138}
\psline[linecolor=black, linewidth=0.04](8.616,-0.42666668)(9.336,0.53333336)
\psellipse[linecolor=black, linewidth=0.04, dimen=outer](11.16,0.47333333)(3.024,1.62)
\psline[linecolor=black, linewidth=0.04](8.136,0.47333333)(14.184,0.47333333)
\psline[linecolor=black, linewidth=0.04](11.16,0.47333333)(10.08,-1.0386667)
\psline[linecolor=black, linewidth=0.04](11.16,0.47333333)(12.024,-1.0386667)
\psframe[linecolor=black, linewidth=0.04, fillstyle=solid,fillcolor=colour0, dimen=outer](12.252,-0.9306667)(11.916,-1.2666667)
\pscircle[linecolor=black, linewidth=0.04, fillstyle=solid,fillcolor=red, dimen=outer](11.13,0.44333333){0.138}
\rput{-8.973528}(0.0889594,2.0203514){\pscircle[linecolor=black, linewidth=0.04, fillstyle=solid,fillcolor=red, dimen=outer](12.918,0.44333333){0.138}}
\pscircle[linecolor=black, linewidth=0.04, fillstyle=solid,fillcolor=red, dimen=outer](9.294,0.50333333){0.138}
\rput[bl](11.076,-0.006666667){$v$}
\pscircle[linecolor=black, linewidth=0.04, fillstyle=solid,fillcolor=red, dimen=outer](14.178,0.44333333){0.138}
\rput{-3.0637553}(-0.01209168,0.43451753){\pscircle[linecolor=black, linewidth=0.04, fillstyle=solid,fillcolor=red, dimen=outer](8.118,0.44333333){0.138}}
\pscircle[linecolor=black, linewidth=0.04, fillstyle=solid,fillcolor=red, dimen=outer](8.598,-0.39666668){0.138}
\psline[linecolor=black, linewidth=0.04](13.416,1.6133333)(12.936,0.53333336)
\pscircle[linecolor=black, linewidth=0.04, fillstyle=solid,fillcolor=red, dimen=outer](13.398,1.5233333){0.138}
\psframe[linecolor=black, linewidth=0.04, fillstyle=solid,fillcolor=colour0, dimen=outer](10.188,-0.9306667)(9.852,-1.2666667)
\psframe[linecolor=black, linewidth=0.04, fillstyle=solid,fillcolor=colour0, dimen=outer](10.248,2.1293333)(9.912,1.7933333)
\psframe[linecolor=black, linewidth=0.04, fillstyle=solid,fillcolor=colour0, dimen=outer](12.288,2.1893334)(11.952,1.8533334)
\rput[bl](2.5426667,-2.2933333){a) $\Delta=5$}
\rput[bl](10.4626665,-2.2933333){b) $\Delta=6$}
\end{pspicture}
}
\caption{Two Halin graphs with maximum size of $LN_d(v)$ for  $\Delta=5$ and $\Delta=6$.}
\label{fig2:LND}
\end{center}
\end{figure}
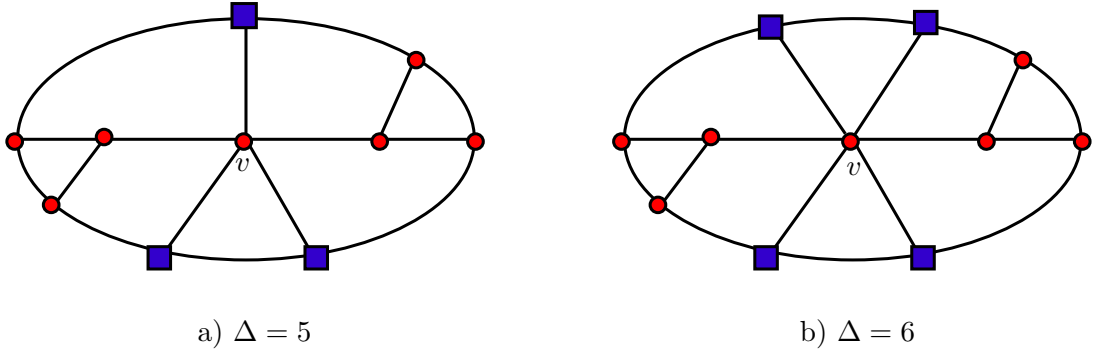

%Let $v$ is a vertex in $T$, and  $v_0,\ldots,v_{\deg(v)-1}$ is an ordering for the neighbours of $v$. We say this ordering is a {\it well-ordering} on the neighbours of $v$ if we have the following conditions.
%\begin{enumerate}
%\item
%If $v\neq r$, then $v_0$ is the parent of $v$.
%\item
%For each $1\leq i<\deg(v)-1$, $\deg(v_i)\leq \deg(v_{i+1})$.
%\item
%If $|E_d(v)|=n_1$ and  $|E_{dc}(v)|=n_2$, then $\{v_1,\ldots,v_{n_1}\}\subseteq E_d(v)$ and $\{v_{n_1-n_2+1},\ldots,v_{n_1}\}\subseteq E_{dc}(v)$.
%\end{enumerate} 
\begin{theorem}
If $G=T\cup C$ is a Halin graph with maximum degree $\Delta$, then 
\[\chi^\prime_s(G)\leq \left\lfloor\frac{3\Delta}{2}\right\rfloor+2.\]
\end{theorem}
\begin{proof}
{
By Theorem~\ref{Halin-cubic}, if $\Delta=3$, then  $\chi^\prime_s(G)\leq \left\lfloor\frac{3\Delta}{2}\right\rfloor+2= 6$.
We now assume that $\Delta\geq 4$ and find a star edge coloring of $G$, as follows.  
We first color some special edges of $C$, then give a star edge coloring of $T$ with at most $\left\lfloor\frac{3\Delta}{2}\right\rfloor$ colors, and finally complete the partial edge coloring of $C$.
 
The partial edge coloring of $C=e_1,\ldots,e_n$ uses colors $\{a,b\}$, with the following patterns. 
Let  $k=n\pmod{3}\times 4$, and use repetative coloring patterns $(a,b,a,\star)$ and $(a,b,\star)$  for edges of paths $e_1,\ldots, e_{k}$ and $e_{k+1},\ldots,e_n$, respectively. 
In these patterns, the colors of edges that are shown with the notation $\star$ is determined after coloring $T$.
For every vertex $v$ of $T$, we denote the number of uncolored edges of $C_d$ that are incident  $v$, with $UC_d(v)$.
Note that if $|LN_d(v)|= \left\lfloor\frac{\Delta}{2}\right\rfloor+1$, then there are two leaves $v_1$ and $v_2$ adjacent to $v$ such that $C$ contains path $v_0v_1v_2v_3$, where edges  $v_0v_1$ and $v_2v_3$ belong to $C_d$. 
Considering the partial edge coloring of $C$, it is clear that $v_0v_1$ or $v_2v_3$ are colored before with a color in $\{a,b\}$; therefore,  $UC_d(v_1)=0$ or $UC_d(v_2)=0$. 
Therefore, there exist at most $\left\lfloor\frac{\Delta}{2}\right\rfloor$ uncolored edges of $C_d$ incident to a vertex   in $LN_d(v)$.
%Thus, we define $E_{dc}(v)\subseteq E_d(v)$ as the set containing the edges with color $a$ or $b$ in all paths   of length 3  with the end edges in $E_d(v)$. 
 
We now give a star edge coloring of $T$ with the color set $\mathcal{C}=\{1,\ldots,\left\lfloor\frac{3\Delta}{2}\right\rfloor\}$. For this purpose, we first choose an arbitrary vertex $r\in V(T)$ of degree $\Delta$ as the root. Let $\ell=0$, and $r$ is a vertex in level 0 of $T$. 
Suppose that for each vertex $v\neq r$ of $T$, $v_0,\ldots,v_{\deg(v)-1}$ are the neighbors of $v$, such that $v_0$ is the parent of  $v$, and  for every $i$, $1\leq i< \deg(v)-1$, $UC_d(v_i)\geq UC_d(v_{i+1})$.
Now we do the following steps to provide a star edge coloring of $T$.

\noindent{\bf  Algorithm 1.} Coloring  characteristic tree $T$ of Halin graph $G=T\cup C$.
\begin{enumerate}
\item[Step 1.] Properly color the incident edges to $r$ with colors $1,2,\ldots,\Delta$, and let $\ell=1$. 
\item[Step 2.] For each vertex $v$ at level $\ell$ which has at least an uncolored incident edge  do
\begin{enumerate}
\item[Step 2.1] Let $m=\min\{\deg(v)-1, \left\lfloor\frac{\Delta}{2}\right\rfloor\}$, and properly color edges  of $\{vv_i:i=1,\ldots,m\}$ with the colors that were not used in the coloring of the edges incident to $v_0$.
\item[Step 2.2]  Properly color the remaining edges incident to $v$ by using any color from $\mathcal{C}$ not yet used in step 7
such that any path of length four in $T$ is not bi-colored.
\end{enumerate}
\item[Step 3.] Stop if all edges are already colored. Otherwise, Let $\ell=\ell+1$ and  go to step 2.
\end{enumerate}

We now ready to complete the partial edge coloring of cycle $C$. 
Let $e=uv$ is an uncolored edges in $C$.
% If $e\in C_s$, then we choose a color in $\mathcal{C}\setminus \Phi(p(v))$. 
If $e\in C_s$, then we   choose a color of $\mathcal{C}$ which is not used for the edges incident to $u_0$. 
 Hence, suppose that $e\in C_d$. Note that if $e=uv$ is an uncolored edge in $C_d$,  then   $UC_d(v),UC_d(u)>0$. Thus,  in the provided star edge coloring of $T$,  the color of  $uu_0$ (resp. $vv_0$) is not used for  any edges incident to $p(u_0)$ (resp. $p(v_0)$) if there exists. 
Moreover,  the color of edge $uu_0$  (resp.  $vv_0$) is used for at most $\left\lceil\frac{\Delta}{2}\right\rceil-1$ incident edges to the neighbors of $u_0$ (resp. $v_0$). If  $\{u_{i_1},\ldots,u_{i_k}\}$  (rasp. $\{v_{i_1},\ldots,v_{i_\ell}\}$)  is the set of neighbors of $u_0$ (resp. $v_0$) that there is an incident edge to every one with the color of $uu_0$ (resp. color of $vv_0$), then the colors of edges in  $E^\prime=\{uu_0, u_{i_1}u_0,\ldots,u_{i_k}u_0\}\cup \{vv_0,v_{i_1}v_0,\ldots,v_{i_\ell}v_0\}$ are forbidden for $e$. Let $\mathcal{C}^\prime$ be the set of colors used for the edges in $E^\prime\cup \{uu_0,vv_0\}$. Thus, the number of the allowed colors for $e$ is 
  \begin{equation}\label{enq1}
  |\mathcal{C}\setminus \mathcal{C}^\prime|\geq\left\lfloor\frac{3\Delta}{2}\right\rfloor-(2+k+\ell)\geq \left\lfloor\frac{3\Delta}{2}\right\rfloor -2\times\left\lceil\frac{\Delta}{2}\right\rceil.
  \end{equation}
 Since $\Delta\geq 4$, $  |\mathcal{C}\setminus \mathcal{C}^\prime|\geq 1$ and  we can choose a color for $e$. 
}
\end{proof}

%%%%%%%%%%%%%%%%%%%%%%%%%%%%%%%%%%%%%%%%%%%%%%%%%%%%%%%%%%%%%%%%%%%
%\bibliographystyle{abbrv}
%\bibliography{bibref}
\setlength{\baselineskip}{0.73\baselineskip}

\end{document}